\documentclass[conference]{IEEEtran}
\IEEEoverridecommandlockouts
\usepackage[utf8]{inputenc}

\usepackage{amssymb,amsmath,amsthm,geometry,wasysym,framed,adjustbox,newclude,cite,graphicx,xcolor}
\usepackage{algorithm}
\usepackage{algpseudocode}
\usepackage{caption}
\usepackage{subcaption}
\theoremstyle{theorem}
\newtheorem{theorem}{Theorem}

\title{PrimeTime: A Finite-Time Consensus Protocol for Open Networks\\
\thanks{This work was supported in part by the National Science Foundation (NSF) under
Grant ECCS-2030251 and CMMI-2024774.}
}

\author{Henry W. Abrahamson$^{\dagger}$ and Ermin Wei$^{\ddagger}$ \thanks{$\dagger$ email address: henryabrahamson2022@u.northwestern.edu.} \thanks{$\ddagger$ email address: ermin.wei@northwestern.edu.} \\
\IEEEauthorblockA{\textit{Department of Electrical and Computer Engineering} \\
\textit{Northwestern University}\\
Evanston, USA} }
\date{August 2022}

\begin{document}

\maketitle

\begin{abstract}
In distributed problems where consensus between agents is required but average consensus is not desired, it can be necessary for each agent to know not only the data of each other agent in the network, but also the origin of each piece of data before consensus can be reached. However, transmitting large tables of data with IDs can cause the size of an agent's message to increase dramatically, while truncating down to fewer pieces of data to keep the message size small can lead to problems with the speed of achieving consensus. Also, many existing consensus protocols are not robust against agents leaving and entering the network. We introduce PrimeTime, a novel communication protocol that exploits the properties of prime numbers to quickly and efficiently share small integer data across an open network. For sufficiently small networks or small integer data, we show that messages formed by PrimeTime require fewer bits than messages formed by simply tabularizing the data and IDs to be transmitted. 
\end{abstract}

\section{Introduction}

In many distributed systems, it is necessary for all the agents in the system to agree on some parameter. Some examples include distributed formation control \cite{listmann2009formation}, optimal routing \cite{madan2006routing}, and distributed kalman filtering for state estimation \cite{carli2008kalman}. This is known as the consensus problem: for some state $x_i$ and agents $i = 1...N$, we wish to drive the system such that $x_i = x_j$ for all $i,j$. \cite{jadb2003consensus} provides the first major theoretical exploration of the average consensus problem in a fully distributed setting. In the years since, many average consensus protocols have been developed, such as \cite{Li2011}, \cite{Gao2012}, as well as the classes of algorithms discussed in \cite{consensustutorial}. These protocols all exhibit asymptotic convergence to the average, i.e. $x_i \to \frac{1}{N} \sum x_j$ as time goes to $\infty$, while only requiring agents to communicate with their immediate neighbors to update their estimates.

However, in systems where quick reaction times may be necessary for safety (such as for collision avoidance for self-driving cars), having a finite convergence time may be preferable over asymptotic convergence. There are protocols for consensus with finite termination, such as \cite{Sundaram2007} and \cite{Sandry2014}. However, \cite{Sundaram2007} involves iteratively calculating out the coefficients of the minimal polynomial of the consensus matrix, while \cite{Sandry2014} requires its gains to be set in terms of the eigenvalues of the graph's weight matrix. Both of these approaches would therefore encounter issues if the graph is non-static, since these values would be changing over time.

Others, such as the protocols described in \cite{Chen2011} and \cite{Liu2022}, rely on a leader-follower style of consensus, which may be unsuited to fully distributed problems in which the leader node might exit the network. There are other finite time protocols of the type described in \cite{Wang2010} and \cite{G_mez_Guti_rrez_2018}, which can cope with dynamic graph topologies. However, these protocols all rely on continuous-time dynamics, and so may not always be feasible to approximate with a discrete time system.

Furthermore, none of the finite time algorithms presented above work on \textit{open} networks, i.e. networks which nodes can freely enter or exit. There are some algorithms that can handle open networks, such as the algorithms in \cite{consensustutorial} that are robust to initial conditions, as well as the first order optimization algorithms presented in \cite{freeman2021healing}, but these methods do not have finite termination. 

Additionally, for some applications, such as intersection management, although it is necessary for the agents to achieve consensus (in this case, on which car(s) can enter the intersection next and what direction(s) they can turn), finding the average or the minimum/maximum value is not necessarily helpful for achieving a meaningful consensus. For example, if each potential action is mapped to some integer, then reaching agreement by computing the average is completely devoid of meaning if the average is not itself an integer. Indeed, computing a single number in general may not fully encapsulate what all of the agents in the system want to do. For our intersection management example, it is obvious that there must be some form of ID-tagging - the system should not only agree that some car will turn left next, but that car $i$ in particular will do so. 

We propose PrimeTime, a prime-number-based finite time consensus protocol with finite termination for integer data that intrinsically includes ID-tags for data. PrimeTime is capable of achieving consensus even on open graphs. Furthermore, PrimeTime allows for an arbitrary desired consensus, not necessarily the average. 

The rest of the paper is organized as follows. First, we present the problem that we are trying to solve, and discuss how it is related to the consensus problem. We then present two versions of PrimeTime, and discuss how they evolve through a simple example in section III. In section IV, we provide some intuition for the consensus speed of PrimeTime, as well as its scalability. Finally, in section V, we present a brief simulation study to show how PrimeTime compares with an equivalent algorithm that does not use prime numbers or prime factorization. 

\section{Formulation}

Suppose that we have a connected undirected graph $\mathcal{G} = \{\mathcal{V,E}\}$, where $\mathcal{V} = \{1,2,\dots N\}$ is the set of nodes in the graph, representing the agents in the system, and $\mathcal{E}$ is the set of edges, denoting the lines of communication. We will use the terms ``agent" and ``node" interchangeably in this work. Let $\mathcal{N}_i$ denote the set of neighbors of agent $i$, and let $\mathcal{N}_i^{(m)}$ denote the $m$-hop neighbors of $i$ (note that $\mathcal{N}_i = \mathcal{N}_i^{(1)}$, and that $\mathcal{N}_i^{(0)} = \{i\}$). 

We define the \textit{inclusive} $m$-hop neighbors of $i$ as $\mathcal{N}_i^{(m)+} = \mathcal{N}_i^{(0)} \cup \mathcal{N}_i^{(1)} \cup \dots \cup \mathcal{N}_i^{(m)}$, that is, the set of agents that can be reached with $m$ hops or fewer. We define the \textit{exclusive} $m$-hop neighbors of $i$ as $\mathcal{N}_i^{(m)-} = \mathcal{N}_i^{(m)} - \big( \mathcal{N}_i^{(0)} \cup \mathcal{N}_i^{(1)} \cup \dots \cup \mathcal{N}_i^{(m-1)} \big) = \mathcal{N}_i^{(m)} - \mathcal{N}_i^{(m-1)+}$, that is, the set of agents that can be reached with $m$ hops and no fewer.  Here $-$ is meant in the set-theoretic sense.

In order to achieve consensus, each agent $i$ in the system wants to build a table $\mathcal{T}_i(k)$ consisting of ordered pairs $(x_j,p_j)$ for all $j$ in the network. $x_j$ is the data over which the system wants to achieve consensus, $p_j$ is a unique identifier for agent $j$, and $k$ is the time index. We assume that each $x_j$ is a strictly positive integer, and that they are all bounded above by some integer $M \geq x_j \forall j$. For example, for basic 4-way intersection management, $M = 3$, with $x_i$ assigned to 1, 2, or 3 if agent $i$ wants to turn right, turn left, or go straight respectively. The algorithms we present are for the case when $x_i$ is a scalar, but can be easily extended for vector $x_i$ by transmitting a vector of messages, with one element for each element of $x_i$.

Once $\mathcal{T}_i(k)$ is complete for all $i$, every agent will know the $x_j$ for all agents in the network. This way of looking at the problem makes it seem more akin to a data sharing or data flooding problem. However, if every $\mathcal{T}_i(k)$ is identical, and each agent has an identical decision-making protocol based on $\mathcal{T}_i(k)$, then the system will have achieved consensus. Note that the system does \textit{not} update any $x_i$ in order to have $x_i \to $ some $\bar{x}$. From that perspective, the system is not necessarily achieving consensus on the data directly, but rather on another variable that depends on $\mathcal{T}_i(k)$ (which in turn depends on the data). 

The advantage of performing consensus by constructing a table like above is that the actual consensus value can have an arbitrary relationship to $\mathcal{T}_i(k)$. It does not need to be an average, maximal or minimal value, or any sort of linear function of the data. In cases with small $N$ and $M$, it may even just be a lookup table. 

\section{Algorithms}

\subsection{PrimeTime}

For simplicity, we will begin by assuming a static graph. In PrimeTime, $p_j$ are set as globally unique prime numbers, so that they can function as identifiers. Initialization must therefore be done in a centralized way to avoid double-assignments of a particular prime, either before the system is deployed, or by designating a leader to assign primes once the system is launched. Adding or removing primes from the system to handle nodes entering or leaving can be done in a distributed way, discussed in section \ref{section_opengraphs}.

In order to complete $\mathcal{T}_i$, at every time step each agent transmits the message 

\begin{equation}
\label{ptmsg}
    m_i(k) = \prod_{\{j : (x_j,p_j) \in \mathcal{T}_i(k)\}}p_j^{x_j},
\end{equation}
i.e., the product of all the primes in agent $i$'s table at time $k$, raised to their associated data's power. Now let $\mathcal{L}_j(k)$ be the set of all agents whose data is included in $m_j(k)$. When messages are received, agents compute a prime factorization to recover the data $x_{l \in \mathcal{L}_j}$, which is stored in the exponent of that data's associated prime. Any new data-prime pairs are added to $\mathcal{T}_i(k+1)$ for the next time step.

\begin{figure}
    \centering
    \includegraphics[width=0.45\textwidth]{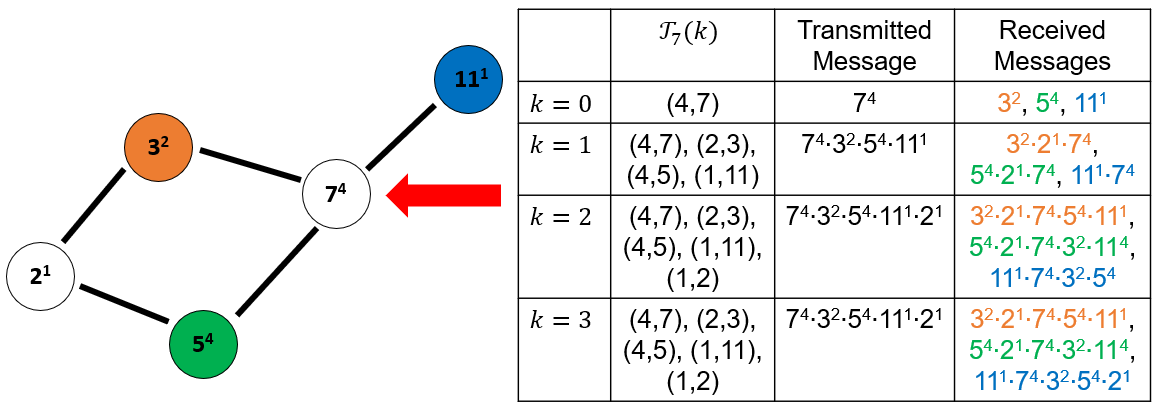}
    \caption{}
    \label{ex1}
\end{figure}

In this way, PrimeTime allows for the encoding of multiple pieces of data, all implicitly ID-tagged, within a single integer. An algorithmic representation of PrimeTime is shown in Algorithm \ref{PrimeTime}, while a brief example of PrimeTime running on a small graph is shown in Figure \ref{ex1}.

\begin{algorithm}
    \caption{PrimeTime}
    \label{PrimeTime}
\begin{algorithmic}
    \State Initialize $\mathcal{T}_i(0) = \{(x_i,p_i)\}$
    \For{$k \geq 0$}
    \State $m_i(k) = \prod_{\{j : (x_j,p_j) \in \mathcal{T}_i(k)\}}p_j^{x_j} $
    \State Transmit $m_i(k)$ 
    \For{$j \in \mathcal{N}_i$}
\State Receive $m_j(k)$
\State Compute the prime factorization of $m_j(k)$ to recover $(x_l,p_l)$ for all $l \in \mathcal{L}_j(k)$ 
\EndFor

\State $\mathcal{T}_i(k+1) = \mathcal{T}_i(k) \cup \{(x_l, p_l) : l \in \bigcup_{j \in \mathcal{N}_i} \mathcal{L}_j(k)\}$

    \EndFor
\end{algorithmic}
\end{algorithm}

The left half of Figure \ref{ex1} shows the graph in question. Each node is labelled with its associated prime, and its data is indicated by the exponent of the prime (e.g. $x_5 = 4$). The right half shows the evolution of PrimeTime from the perspective of node 7, indicated by the red arrow. The leftmost column indicates node 7's local table $\mathcal{T}_7(k)$, while the two columns to the right indicate node 7's transmitted message, $m_7(k)$, and node 7's incoming messages, $m_j(k)$ with $j \in \mathcal{N}_7$. Incoming messages are color coded with their node of origin; e.g., the messages that node 7 receives from node 5 are indicated in green. 

Note that, in steady state, since $\mathcal{T}_i$ will contain the data of every agent in the network, each agent's message will be the product of every agent's prime raised to the corresponding data's power - a potentially large message with large amounts of redundancy. 

\subsection{Incremental PrimeTime}

If dropped packets or other communication errors are common, and if both the range of the data and the network itself are quite small, this redundancy may be welcome. However, if either the data or the network is large, then it might be infeasible or impractical to constantly transmit possibly large integers. In this case, we introduce Incremental PrimeTime, a modified version of PrimeTime that has better scalability through reduced redundancy, shown in Algorithm \ref{incPrimeTime}.

\begin{algorithm}
    \caption{Incremental PrimeTime}
    \label{incPrimeTime}
\begin{algorithmic}
    \State Initialize $\mathcal{T}_i(0) = \{(x_i,p_i)\}, \mathcal{T}_i(-1) = \emptyset$
    \For{$k \geq 0$}
    \State Set $\mathcal{A}_i(k) = \mathcal{T}_i(k) - \mathcal{T}_i(k-1)$
    \State $m_i(k) = \prod_{\{j : (x_j,p_j) \in \mathcal{A}_i(k)\}}p_j^{x_j} $
    \State Transmit $m_i(k)$ 
    \For{$j \in \mathcal{N}_i$}
\State Receive $m_j(k)$
\State Compute the prime factorization of $m_j(k)$ to recover $(x_l,p_l)$ for all $l \in \mathcal{L}_j(k)$ 
\EndFor

\State $\mathcal{T}_i(k+1) = \mathcal{T}_i(k) \cup \{(x_l, p_l) : l \in \bigcup_{j \in \mathcal{N}_i} \mathcal{L}_j(k)\}$
    \EndFor
\end{algorithmic}
\end{algorithm}

In Incremental PrimeTime, instead of transmitting the product of their entire table, agents only transmit new data. More precisely, instead of setting $m_i(k)$ according to (\ref{ptmsg}), each agent first constructs an auxiliary set

\begin{equation}
    \label{auxsetdef}
    \mathcal{A}_i(k) = \mathcal{T}_i(k) - \mathcal{T}_i(k-1).
\end{equation}

Then, $m_i(k)$ is formed by 

\begin{equation}
    \label{iptmsg}
    m_i(k) = \prod_{\{j : (x_j,p_j) \in \mathcal{A}_i(k)\}}p_j^{x_j}.
\end{equation}

\begin{figure}
    \centering
    \includegraphics[width=0.45\textwidth]{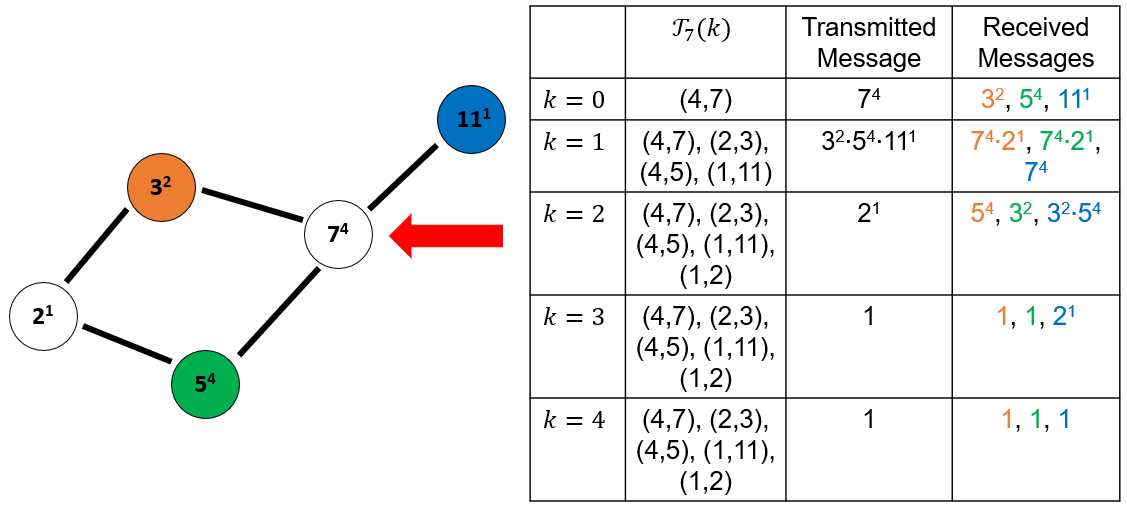}
    \caption{}
    \label{ex2}
\end{figure}

To demonstrate, the same example from Figure \ref{ex1} is shown in Figure \ref{ex2}, but with the system running Incremental PrimeTime instead. Compared to the previous example, node 7 updates its table identically. However, by using $\mathcal{A}_i(k)$ to form messages, each agent only transmits each data-prime pair once, during the time step right after that agent first receives it. This means that, for Incremental PrimeTime, in steady state $m_i(k) = 1$. Once every agent knows every other agent's data, no new data will be obtained, so $\mathcal{A}_i(k) = \emptyset$. 

Incremental PrimeTime has a clear advantage over PrimeTime in that its messages will be smaller, which will allow it to scale better for larger networks and larger $M$. However, PrimeTime has the benefit of having highly redundant messages. Each agent transmits its entire table at each time step, compared to Incremental PrimeTime, in which agents transmit each data-prime pair exactly once. Intuitively, this means that if the system suffers from packet loss or something else that causes the graph topology to be dynamic, PrimeTime will be able to pass along data and complete its tables more consistently than its incremental version.

\subsection{Open Graphs}\label{section_opengraphs}

PrimeTime and Incremental PrimeTime can both be easily extended to handling open graphs, assuming that the system is in steady state. To add a new agent $i$ to the network, $i$ just needs to query one of its neighbors for its table to find the smallest unused prime. It then starts performing PrimeTime as if it were at time $k = 0$, and so all other agents will eventually update their table with the new prime-data pair.

Now, consider the case in which node $i$ wants to leave the network. Since the data is upper bounded by $M$, we can set $x_i = M'$, the smallest integer that is not used in the range of $x_i$, as an indicator of leaving the network. Just before agent $i$ leaves, it transmits $m_i(k)' = m_i(k) \times p_i^{M'}$, where $m_i(k)$ is formed in the usual way according to PrimeTime or Incremental PrimeTime. When agent $j$ receives this message and computes the prime factorization to recover $p_i^{M'}$, it simply removes the pair $(x_i,p_i)$ from its table and includes $p_i^{M'}$ in its product for its next message $m_j(k+1)$, to pass agent $i$'s ``goodbye" along the network. 

These methods are only guaranteed to work for addition and removal of agents once the system has achieved steady state. If the system has not yet achieved steady state, then problems may arise if the graph's topology contains any loops.

\section{Performance Analysis}

\subsection{Consensus Speed}

We start by rewriting the equation for PrimeTime messages (\ref{PrimeTime}) in terms of its graph theoretic representation. At time $k=0$, agents transmit only their own prime and data, and so receive the data from their one-hop neighbors. At time $k=1$, agents transmit both the primes and data of their 1-hop neighbors and their own, and so receive that data from their 1-hop neighbors' 1-hop neighbors, i.e. their two hop neighbors. As such, we can see that agents hear back from their $k$-hop neighbors at time $k-1$, and we can write the following equation for $\mathcal{T}_i(k)$:

\begin{equation}
    \mathcal{T}_i(k) = \bigcup_{j \in \mathcal{N}_i^{(k)+}}(p_j,x_j),
    \label{N_table_pt}
\end{equation}

from which we can derive the equation for $m_i(k)$,

\begin{equation}
    m_i(k) = \prod_{j \in \mathcal{N}_i^{(k)+}}p_j^{x_j}.
    \label{N_msg_pt}
\end{equation}

We now present the finite termination time for PrimeTime.

\begin{theorem}
Let d be the diameter of $\mathcal{G}$. Then for all i, $\mathcal{T}_i(k)$ contains the data for all agents in the network for all $k \geq d$, and no sooner. 
\end{theorem}

\begin{proof}
First, we will show that every agent's table is complete at time $d$. Let $k \geq d$, and consider some node $i$. By definition of the diameter of a graph, $d$ is the smallest integer such that all other nodes can be reached from node $i$ in $d$ hops or fewer. Therefore, $\mathcal{N}_i^{(d)+} = \mathcal{V}$, so $\mathcal{N}_i^{(k)+} = \mathcal{V}$. Plugging this into \eqref{N_table_pt} yields $\mathcal{T}_i(k) = \bigcup_{j \in \mathcal{N}_i^{(k)+}}(p_j,x_j) = \bigcup_{j \in \mathcal{V}}(p_j,x_j)$, which means that the table is complete. 

Now, let $k < d$. Again, by definition of the diameter, there exists some pair of nodes $i,j$ such that $j$ is a $d$-hop neighbor of $i$, and $j$ is not a $\ell$-hop neighbor of $i$ for any $\ell < d$. Therefore, $(x_j,p_j) \notin \mathcal{T}_i(k)$, and so not all tables are complete.
\end{proof}

Immediately from Theorem 1, we see that $m_i(k) = \prod_{i \in \mathcal{V}}p_i^{x_i}$ for all $k \geq d$. This means that every agent will transmit all agents' data from time $d$ onwards. If the network is large, then this message may be impractical to form and compute the prime factorization of. 

Now, consider the messages for Incremental PrimeTime. In this case, messages are only constructed from new entries into $\mathcal{T}_i(k)$. However, the same rate of information propagation as before holds, and so $\mathcal{T}_i(k)$ updates as in (\ref{N_msg_pt}). This means that Incremental PrimeTime also completes its table in exactly $d$ time steps, using an identical argument as above. The messages, however, update a bit differently:

\begin{equation}
    m_i(k) = \prod_{j \in \mathcal{N}_i^{(k)-}}p_j^{x_j}.
    \label{N_msg_ipt}
\end{equation}

Because $\mathcal{N}_i^{(d)+} = \mathcal{V}$ for all $i$, $\mathcal{N}_i^{(d+1)-} = \mathcal{N}_i^{(d+1)}-\mathcal{N}_i^{(d)+} = \emptyset$. Therefore, all agents will transmit just $1$ at time $d+1$ and above. As before, tables will be complete after $d$ rounds of communication, after which consensus can be achieved.

\subsection{Scalability}

Scalability appears as an immediate concern for PrimeTime, since messages that are composed of potentially long products might get large. This might require custom implementations of large, unsigned integers, and may also make the prime factorization step infeasible, depending on hardware constraints of the application of interest. Certainly, Incremental PrimeTime would scale better than PrimeTime in this regard. We investigate the scalability of Incremental PrimeTime in simulation in Section V, but provide a few intuitions below.

If the number of agents in the network $N$ is large, then agents will be forced to use larger and larger prime numbers. The steady-state message size for PrimeTime will therefore increase at least as fast as $N$ factorial, if not worse. A large range of integer data would also cause the message size to blow up quickly, since some of the primes would be raised to large powers. Although Incremental PrimeTime doesn't have the same steady-state message as PrimeTime, messages in the transient could still be quite large, for the same reasons as above. 

The topology of the graph will also affect the message size - in the case of PrimeTime, it affects how fast the messages increase in size, while for Incremental PrimeTime, it has a direct connection the transient message size. Recall that Incremental PrimeTime messages can be written according to \eqref{N_msg_ipt}. Therefore, the message size for agent $i$ at time $k$ is directly related to $|\mathcal{N}_i^{(k)-}|$, the number of exclusive $k$-hop in-neighbors of $i$. All else being equal, graphs with large, highly connected clusters would therefore have the largest message sizes under Incremental PrimeTime. 

\section{Numerical Results}

\begin{table*}
    \centering
    \begin{tabular}{|c|c|c|c|} \hline
       Graph Parameters & \begin{tabular}{@{}c@{}}Max Bytes \\ (Incremental PrimeTime)\end{tabular} & \begin{tabular}{@{}c@{}} Avg Bytes \\ (Incremental PrimeTime) \end{tabular} & \begin{tabular}{@{}c@{}}Avg Bytes \\ (Vectorized Scheme)\end{tabular} \\  \hline
       $N = 15$, $r = 0.36$, $M = 3$ & 14 & 5.44 & 11.54\\ \hline
       $N = 15$, $r = 0.36$, $M = 5$ & 19 & 6.94 & 11.71 \\ \hline
       $N = 10$, $r = 0.36$, $M = 5$ & 12 & 4.98 & 8.53 \\ \hline
       $N = 15$, $r = 0.5$, $M = 3$ & 16 & 6.88 & 16.67 \\ \hline
    \end{tabular}
    \caption{Average byte requirements for messages formed by Incremental PrimeTime and the vectorized scheme, along with the maximum byte requirement for Incremental PrimeTime for reference.}
    \label{ipt_table}
\end{table*}

To investigate the possible scalability of Incremental PrimeTime, we performed a brief simulation study. We leave out basic PrimeTime, as Incremental PrimeTime achieves the same convergence speed with less communication overhead. For each run, we generated 100 random geometric graphs, since graphs of this type are generally a good model for many physical applications (e.g. in environmental sensor networks, physical distance is often the main determiner of whether two nodes can communicate or not \cite{carli2008kalman}). $N$ points were uniformly generated on the unit square, and edges between two nodes were added if they were within $r$ distance of each other. If a graph was not connected, we regenerated it with the same parameters. We then compiled the non-1 messages formed by Incremental PrimeTime over all the graphs, and saw which messages fit within 4-byte unsigned integers versus 8-byte unsigned integers, and used that to compute the average amount of bytes used per message. We chose 8-bytes as a cutoff point, as the largest data type natively implemented for many common programming languages is the 64-bit unsigned integer, and because from that point on prime factorizations start to become more costly computationally. 

As a point of comparison, we also ran another algorithm on these graphs - it still builds up tables $\mathcal{T}_i(k)$ with data and IDs, but instead of using prime numbers and prime factorizations, it simply transmits the data-ID pairs as two integers, and transmits multiple pairs as a long vector (e.g., 2 pairs would require four integers, 3 pairs would require six, etc.). This vectorized scheme mimics the behavior of Incremental PrimeTime exactly, except for the message size, allowing us to see how much communication cost we save using our prime factorization scheme. Note that we assumed that the vectorized scheme used 16-bit integers, since that is the default for C and C++ (for reference, an int in Python is typically 32 bits). 

\begin{figure*}
    \centering
    \begin{subfigure}[b]{0.45\textwidth}
        \includegraphics[width=\textwidth]{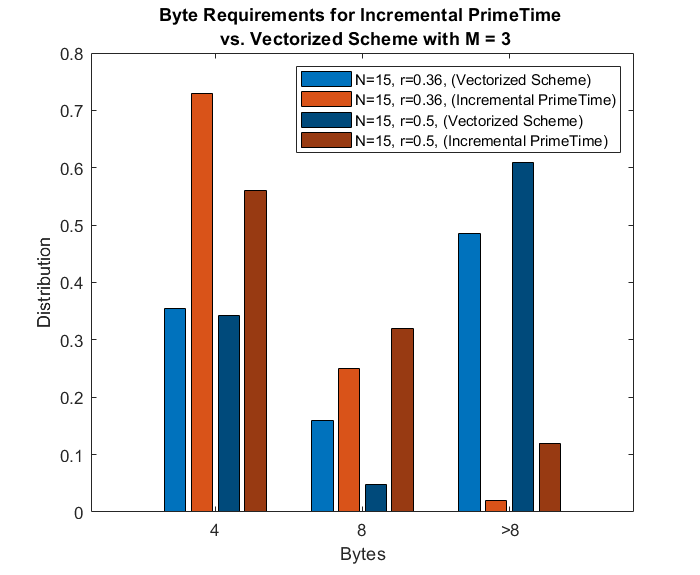}
        \caption{}
        \label{compm3}
    \end{subfigure}
    \begin{subfigure}[b]{0.45\textwidth}
        \includegraphics[width=\textwidth]{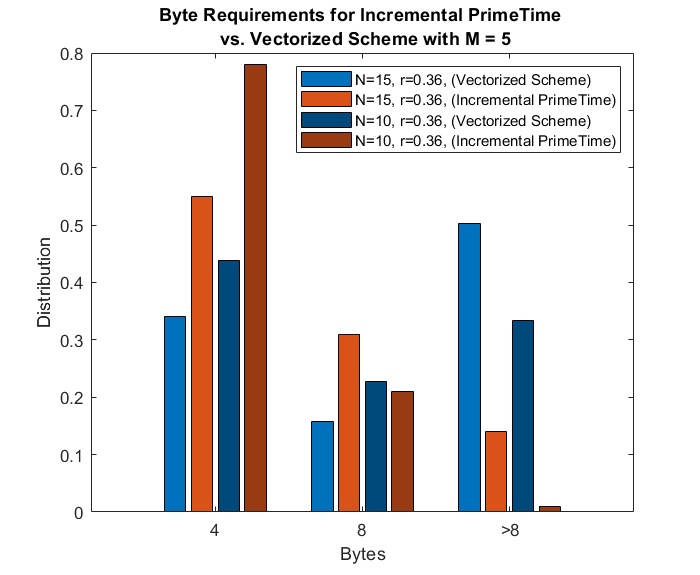}
        \caption{}
        \label{compm5}
    \end{subfigure}
    \caption{Message sizes for Incremental Primetime and the vectorized scheme when $M=3$ (a) and $M=5$ (b).}
    \label{compfig}
\end{figure*}

The results are shown in Figure \ref{compfig} and Table \ref{ipt_table}. Table \ref{ipt_table} shows the average message size for both Incremental PrimeTime and the vectorized scheme. It also includes the bytes needed for the largest message formed by Incremental PrimeTime as a rough measure of practicality. Figure \ref{compfig} shows the distribution of the byte requirements for messages formed by both schemes, with Incremental PrimeTime shown in red and the vectorized scheme shown in blue.

As seen in Figure \ref{compm3}, for random geometric graphs with $N = 15$, $r = 0.36$, and a maximum data value of $M = 3$, $\sim$98\% of messages formed by Incremental PrimeTime were able to fit within an 8-byte unsigned integer or smaller. For contrast, transmitting a single data-ID pair requires 4 bytes, with the average byte requirement for the vectorized scheme in this case being about 2.1 times that of Incremental PrimeTime. When we raise $M$ to 5 (Figure \ref{compm5}), Incremental PrimeTime still maintains a smaller average byte requirement of 6.94 (compared to 11.71), but the proportion of messages that require more than 8 bytes rises to 13\%. The average message size for  the vectorized scheme is about the same between the two previous cases, since the integer data is being transmitted as a 16-bit integer, regardless of its size.

On the other hand, reducing the network size to $N = 10$ while keeping $M=5$ lowers the proportion of messages larger than 8 bytes to $<1\%$, making an 8 byte integer implementation seem more reasonable, while still having messages about 1.7 times smaller than the vectorized scheme. Lastly, to demonstrate how the graph's topology affects the message size, we set $N$ and $M$ back to $15$ and $3$ and set $r=0.5$ (back to Figure \ref{compm3}). In this case, although 12\% of messages were above 8 bytes, the relative reduction in average message size increased to about 2.4 times. 

\section{Conclusion}

We have shown how PrimeTime uses the properties of prime numbers to efficiently encode multiple distinct pieces of implicitly ID-tagged data within a single message, allowing for networks to achieve finite-time consensus while still accommodating open graphs. PrimeTime seems well-suited to applications with small graphs and a limited data range that want global information sharing, such as intersection management for autonomous vehicles. However, PrimeTime does have issues with scalability, and so we introduce Incremental PrimeTime to help lower the message size. Future work may include extending PrimeTime to allow for directed graphs and dynamic graph topologies, as well as looking into further message truncation to keep the messages within 8 bytes for practical implementations.

\section{Acknowledgments}

We would like to thank Anthony Goeckner, Qi Zhu, and Randy Berry for their helpful discussions and insights.

\medskip

\bibliographystyle{ieeetr}
\bibliography{ref}

\end{document}